\documentclass[11pt,twoside]{article}
\usepackage[left=3cm,right=3cm,top=3cm,bottom=3cm]{geometry}
\usepackage{amssymb}
\usepackage{graphicx}
\usepackage{graphicx,type1cm,eso-pic,color}
\usepackage{amssymb}
\usepackage{amssymb,amsmath}
\usepackage{graphicx,epsfig}
\usepackage{enumerate}
\usepackage{color}
\usepackage {multicol}
\usepackage{lineno}
\numberwithin{equation}{section}
\newtheorem{theorem}{Theorem}

\newtheorem{proposition}{Proposition}
\newtheorem{remark}{Remark}

\newtheorem{example}{Example}

\makeatletter
\AddToShipoutPicture{
	\setlength{\@tempdimb}{.5\paperwidth}
	\setlength{\@tempdimc}{.5\paperheight}
	\setlength{\unitlength}{1pt}
	\put(\strip@pt\@tempdimb,\strip@pt\@tempdimc)
}
\makeatother

\usepackage{scalerel,stackengine}
\stackMath
\newcommand\reallywidehat[1]{%
	\savestack{\tmpbox}{\stretchto{%
			\scaleto{%
				\scalerel*[\widthof{\ensuremath{#1}}]{\kern-.6pt\bigwedge\kern-.6pt}%
				{\rule[-\textheight/2]{1ex}{\textheight}}
			}{\textheight}%
		}{0.5ex}}%
	\stackon[1pt]{#1}{\tmpbox}%
}
\begin{document}
	\setcounter{page}{1}
	\thispagestyle{empty}
	\markboth{}{}

	\pagestyle{myheadings}
	\markboth{}{ }
	
	\date{}
	
	
	\noindent  
	
	\vspace{.1in}
	
	{\baselineskip 20truept
		
		\begin{center}
			{\Large {\bf Extropy and Varextropy estimators with applications
			}} \footnote{\noindent	{\bf * } Corresponding author E-mail: skchaudhary1994@kgpian.iitkgp.ac.in\\
				{\bf **} E-mail: nitin.gupta@maths.iitkgp.ac.in
			}			
			
		\end{center}
		
		\vspace{.1in}
		
		\begin{center}
			{\large {\bf  Santosh Kumar Chaudhary*, Nitin Gupta**}}\\
			{\large {\it Department of Mathematics, Indian Institute of Technology Kharagpur, West Bengal 721302, India }}
			\\
		\end{center}
	}
	\vspace{.1in}
	\baselineskip 12truept

	
	\begin{center}
		{\bf \large Abstract}\\
	\end{center}
	In many statistical studies, the measure of uncertainties like entropy, extropy, varentropy and varextropy of a distribution function is of prime interest. This paper proposes estimators of extropy and varextropy. Proposed estimators are consistent. Based on extropy estimator, a test of symmetry is given.  The proposed test has the advantage that we do not need to estimate the centre of symmetry. The critical value and power of the proposed test statistics have been obtained. The test procedure has been implemented on six real-life data sets to verify its performance in identifying the symmetric nature. \\
	\\
	\textbf{Keyword:} Extropy, Estimator,  Symmetry, Varextropy, Weighted Varextropy.  \\
	\\
	\noindent  {\bf Mathematical Subject Classification}: {\it 62B10, 62D05}
	\section{Introduction}\label{s1intro}
	
	Let X be an absolutely continuous random variable with probability density function(pdf) f(x). Let $l_X  = inf\{x \in \mathbb{R}: F(x) > 0\}, u_X  = sup \{x \in \mathbb{R}: F(x) < 1\}$ and $S_X  = (l_X , u_X ).$ Shannon (1948) defined differential entropy as a measure of uncertainty as	
	\begin{align*}
		H(X)=E(-\log(f(X)))=-\int_{S_X} f(x) \log(f(x))dx
	\end{align*}
	The expectation of the information content of an absolutely continuous random variable is represented by Shannon's entropy, a well-known information measure. Information theory applications use the corresponding variance, known as varentropy. Numerous works by Bobkov and Madiman (2011), Kontoyiannis and Verdu (2014), Arikan (2016), Di Crescenzo and Paolillo (2021), and Maadani et al. (2021) have recently made contributions on varentropy.	Varentropy of  $X$
	is defined as 	
	\begin{align*}
		VH(X)&=Var(-\log(f(X)))=Var(\log(f(X)))  \\
		&=E(\log(f(X)))^2- [E(\log(f(X)))]^2		\\
		&=\int_{S_X} f(x) (\log(f(x)))^2dx - \left[\int_{S_X} f(x) \log(f(x))dx\right]^2	
	\end{align*}
	
	This varentropy measure is commonly used in data compression, finite blocklength information theory, and statistics as it aids in determining the ideal code length for data compression, source dispersion, and other relevant considerations. Furthermore, statistics, have proven to be a superior alternative to the kurtosis measure for continuous density functions (see Arikan (2016) and Maadani et al. (2022)).

	An alternative measure of uncertainty, extropy of a non-negative absolutely continuous random variable X defined by Lad et al. (2015) is given as
	\begin{align}
		J(X)=E\left(-\frac{1}{2}f(X)\right)=-\frac{1}{2}\int_{S_X} f^2(x)dx.
	\end{align}
	
	Entropy was first introduced by Shannon (1948) and it measures the average level of uncertainty related to the results	of the random experiment. The Shannon entropy for continuous random variable $X$ is defined  as:
	\begin{equation}
		\label{1eq2}
		H(X) = - {\int_{S_X} {f(x) \ln \left( f(x)\right)dx}}=E\left(- \ln f(X)\right).
	\end{equation} 
	Lad et al. (2015) defined the compliment dual of the Shannon entropy called extropy. The extropy of continuous random variable $X$ is defined as:	
	\begin{equation}\label{extropy}
		J(X)=-\frac{1}{2} \int_{S_X} f^2(x)dx=-\frac{1}{2}E\left(f(X)\right).
	\end{equation}
	The cumulative residual extropy  of continuous random variable $X$ is defined as :
	\begin{align}
		\xi J (X) &=- \frac{1}{2} \int_{S_X} \bar{F}^2 (x) dx.
	\end{align} 
	The cumulative past extropy of $X$ is defined as:
	\begin{align}
		\bar{\xi} J (X) &=- \frac{1}{2} \int_{S_X} F^2 (x) dx.
	\end{align}	
	\indent The concept of $k$-records was introduced by Dziubdziela and Kopociński  (1976); (also see  Ahsanullah (1995) and Arnold et al. (1998)). The pdf of the $n$th upper $k$-record value $U_{n,k}$ and the $n$th lower $k$-record value $L_{n,k},$ respectively, are given by (see Arnold et al. (2008) and  Ahsanullah (2004)) 	
	\begin{align*}
		&f_{U_{n,k}}(x)=\dfrac{k^n}{(n-1)!}[-\log \bar{F}(x)]^{n-1} \bar{F}(x)^{k-1}f(x),\  x\in\mathbb{R}, \label{r1}\\
		\text{and} \ \ \ &f_{L_{n,k}}(x)=\dfrac{k^n}{(n-1)!}[-\log F(x)]^{n-1} F(x)^{k-1}f(x), \  x\in\mathbb{R}. 
	\end{align*}	
	\noindent The cdf of $U_{n,k}$ and  $L_{n,k}$, respectively, are  
	\begin{align*}
		F_{U_{n,k}}(x)&=  1-\bar{F}^k (x)   \sum_{i=0}^{n-1} \frac{(-k \log\bar{F} (x) )^i}{i!},\\
		\text{and} \ \ \ 	F_{L_{n,k}}(x)&=  {F}^k (x)   \sum_{i=0}^{n-1} \frac{(-k \log F (x) )^i}{i!}.
	\end{align*}
	
	Extropy and its applications have recently attracted the attention of various researchers. Entropy and extropy in machine learning are two of the many techniques and concepts that are being used to solve complex problems easily. One advantage of extropy as compared to other measures is that it yields an expression for finite mixture distributions, while no such expression is available for closed-form entropy and variance measures (see Toomaj et al. (2023)). The most important advantage of extropy is that it is easy to compute, and it will therefore be of great interest to explore its important potential applications in developing goodness-of-fit tests and inferential methods. Tahmasebi and Toomaj  (2020) studied the stock market in OECD countries based on a generalization of extropy known as negative cumulative extropy. Balakrishnan et al. (2022) applied another version of extropy known as the Tsallis extropy to a pattern recognition problem. Kazemi et al.  (2021) explored an application of a generalization of extropy known as the fractional Deng extropy to a problem of classification. Tahmasebi et al. (2022) used some extropy measures for the problem of compressive sensing. 
	
	This paper is organized as follows. We discuss about varextropy in Section 2. We define weighted vaextropy in Section 3.  In Section 4, we propose varextropy estimators. We obtained a characterization of uniform distribution using varextropy in Section 5. In Section 6, a test of symmetry is proposed. Section 7 concludes this paper. 
	
	\section{Varextropy}

	Two random variables can have the same extropy in some situations, which prompts the age-old question, "Which of the extropies is a more appropriate criterion for measuring the uncertainty ?". For example, consider random variables $U$ and $V$ (see Balakrishanan et al 2020) with pdf's 
	
	\begin{eqnarray*}
		f_{U}(x)=
		\begin{cases}
			1, \hspace{4mm} 0<x<1\\
			0, \hspace{5mm} otherwise
		\end{cases}~~~~~~~~~~~~~~~~~
		f_{V}(x)=
		\begin{cases}
			2e^{-2x},\hspace{4mm} x>0\\
			0,\hspace{6mm} otherwise
		\end{cases}
	\end{eqnarray*}
	
	We get $J(U)=J(V)=-1/2$ and $VJ(U)=0 $ and $VJ(V)=1/12.$ This is the motivation behind considering the variance of $-\frac{1}{2}f(x),$ which is known as varextropy of a random variable $X.$ So, varextropy can also play a role as a measure of uncertainty. Varextropy of absolutely continuous random variables $X$ is defined as (see  Vaselabadi et al. (2021), Goodarzi (2021) and Zaid et al. (2022))

	\begin{align}
		VJ(X)=Var\left(-\frac{1}{2}f(X)\right) \nonumber &=E\left(-\frac{1}{2}f(X)-J(X)\right)^2 \nonumber\\
		&=\frac{1}{4} E(f^2(X))- J^2(X) \nonumber\\
		&=\frac{1}{4}E(f^2(X))- \frac{1}{4}[E(f(X))]^2 \nonumber\\
		&=\frac{1}{4}  \int_{S_X} f^3(x)dx - \frac{1}{4} \left[ \int_{S_X} f^2(x)dx \right]^2
	\end{align}
	
	Suppose that $X_1,\dots, X_n$ are independent and identically distributed observations with cdf $F$ and pdf $f.$ An observation $X_j$ will be called an upper record value if its value exceeds that of all previous observations. Thus, $X_j$
	is an upper record	if $X_j >X_i$ for every $j > i.$ An observation $X_j$ will be called a lower record value if its value is less than that of all previous observations. Thus, $X_j$
	is an lower record	if $X_j < X_i$ for every $j > i.$ Belzunce et al. (2001) showed that if $X \leq_{disp} Y$, then $U^X_n \leq_{disp} U^Y_n,$ where $U^X_n$
	and $U^Y_n$	are the nth upper records of X and Y, respectively. . Qiu (2017) showed that if $X \leq_{disp} Y$, then $J(X) \leq J(Y)$ and $ J(U_n^X) \leq J(U^Y_n).$ Vaselabadi et al. (2021) showed that if $X \leq_{disp} Y$, then $VJ(X) \geq VJ(X).$ In view of these results, we obtain the following proposition immediately.
	
	The values of varextropy for some standard distribution are given below, for more examples see, Vaselabadi et al. (2021).
	\begin{example}
		When $X$ has a uniform distribution on the interval (0,1), then 	$VJ(X)=0.$ 
	\end{example}
	
	\begin{example}
		Let random variable X have exponential distribution with cdf $F_X(x)=1- e^{-\lambda x}, \ x >0.$ Then 	$VJ(X)=\lambda^2/48, J(X)=\frac{\lambda^2}{16}$ and $VJ(U_n^X)= J(X) \left[\frac{4\Gamma(3n-2)}{\Gamma^3(n) 3^{3n-2}}- \frac{\Gamma^2(2n-1)}{\Gamma^4(n) 4^{2n-2}}\right].$
	\end{example}
	
	\begin{example}
		When $X$ has normal distribution with mean $\mu$ and variance $\sigma^2,$ then 	$VJ(X)=\frac{2-\sqrt{3}}{16\pi\sigma^2 \sqrt{3}}.$ 
	\end{example}

	\begin{proposition}
		Let $X\leq_{disp} Y$ then $VJ(U_n^X) \geq VJ(U_n^Y).$
	\end{proposition}
	
	\begin{proposition}
		RVs X and Y are identically distributed then $VJ(U^X_n)=VJ(U^Y_n).$
	\end{proposition}
	
	\begin{proposition}
		RVs X and Y are identically distributed then $VJ(X_{m:n})=VJ(Y_{m:n}).$
	\end{proposition}

	Note that $VJ(X)\geq 0,$ for any Random variable $X.$ Vaselabadi et al. (2021) obtained several varextropy properties as well as conditional varextropy properties based on order statistics, record values, and proportional hazard rate models. The article contains some comparative results regarding varextropy and varentropy. Goodarzi (2022) provided lower bounds for varextropy, obtained the varextropy of a parallel system, and used the varextropy of order statistics to construct a symmetry test. Zaid (2022) computed the entropy, Varentropy, and Varextropy measures in closed form for generalized and q-generalized extreme value distributions. Varentropy is sometimes independent of the model parameters, whereas the varextropy measure is more adaptable, for example when $X$ has a normal distribution with mean $\mu$ and $\sigma^2$ (see Vaselabadi et al. (2021)).

	The main purpose of this paper is to estimate the varextropy of a continuous random variable due to the interesting properties and potential applications of varextropy and to test the uniformity using the second estimator we proposed. Here, we are estimating the varextropy of a continuous random variable.

	\section{Weighted varextropy}	
	Applications of the weighted distributions include distribution theory, dependability, probability, ecology, biostatistics, and applied statistics. Two random variables can have the same extropy as well as same varextropy in some situations. For example, consider random variables $X$ and $Y$ with pdf's respectively, 
	\begin{eqnarray*}
		f_{X}(x)=
		\begin{cases}
			2x, \hspace{4mm} 0<x<1\\
			0, \hspace{5mm} otherwise
		\end{cases}~~~~~~~~~~~~~~~~~
		f_{Y}(x)=
		\begin{cases}
			2(1-x),\hspace{4mm} 0<x<1\\
			0,\hspace{6mm} otherwise
		\end{cases}
	\end{eqnarray*}
	We get $J(X)=J(Y)=-2/3,$ \  $VJ(X)=VJ(Y)= 1/18$ but $VJ^w(X)=1/12$ and $VJ^w(Y)= 1/180$ So here weighted varextropy can also play a role as a measure of uncertainty. Analogous to Gupta and Chaudhary (2023), we define general weighted varextropy as
	\begin{align}\label{defweightedvarextropy}
		VJ^w(X)&=Var\left(-\frac{1}{2}w(X)f(X)\right) \nonumber\\
		&=\frac{1}{4}\left[E(w^2(X)f^2(X))- (E(w(X)f(X)))^2\right] \nonumber\\
		&=\frac{1}{4} \left[ \int_{S_X}w^2(x)f^3(x)dx - \left(\int_{S_X} w(x)f^2(x)dx\right)^2\right]
	\end{align}
	
	When $w(x)=x,$ then weighted varextropy is given as 
	\begin{align*}
		VJ^x(X)&=\frac{1}{4} \left[ \int_{S_X}x^2f^3(x)dx - \left(\int_{S_X} xf^2(x)dx\right)^2\right]
	\end{align*}
	
	\begin{remark}
		When $w(x)=1,$ then $VJ^w(X)=VJ(X).$
	\end{remark}
	
	\begin{example}
		If X has a uniform distribution on the interval (a,b), then $VJ^x(X)=\frac{1}{48}.$ Note that $VJ^x(X)$ is free of both $a$ and $b.$
	\end{example}
	
	\begin{example}
		Let random variable X have exponential distribution with cdf $F_X(x)=1- e^{-\lambda x}, \ x >0.$ Then $VJ^x(X)=\frac{5}{1728}.$ Note that $VJ^x(X)$ is free of $\lambda.$
	\end{example}

	\section{Varextropy Estimators}
	Here, we'll first discuss a few varextropy estimators. Let $X_1, X_2, X_3, \dots, X_n$ be a random sample from a distribution with unknown probability density function (pdf) f and (cdf) F. Suppose $X_{1:n}, X_{2:n}, X_{3:n}, \dots, X_{n:n}$ order statistics based on random sample $X_1, X_2, X_3, \dots, X_n.$ The empirical distribution function of cdf $F$ is defined as 		
	\begin{eqnarray*}
		{F}_{n}(x)=
		\begin{cases}
			0, \hspace{4mm} x<X_{1:n}\\
			\frac{i}{n}, \hspace{4mm} X_{i:n}\leq x<X_{i+1:n}, \ \ \ i=1,2,\dots, n-1.\\
			1, \hspace{5mm} x \geq X_{n:n}.
		\end{cases}
	\end{eqnarray*}

	\subsection{The first estimator}
	
	The varextropy of X can be written as 
	
	\begin{align*}
		VJ(X)=\frac{1}{4} \int_{0}^{1} \left(\frac{d}{dp}(F^{-1}(p))\right)^{-2} dp- \frac{1}{4} \left[\int_{0}^{1} \left( \frac{d}{dp}(F^{-1}(p))\right)^{-1}dp\right]^2
	\end{align*}
	Following the idea of Vasicek (1976) an estimator of  $VJ(X)$  will be calculated by replacing the distribution function $F$ with empirical distribution function ${F}_n$ and using the difference operator in place of a differential operator. The derivative of $F^{-1} (p)$ with respect to $p$, that is, $\frac{dF^{-1} (p)}{dp}$ will be estimated as 
	\begin{align*}
		\frac{X_{i+m:n}-X_{i-m:n}}{{F}_n(X_{i+m:n})-{F}_n(X_{i-m:n})}=\frac{X_{i+m:n}-X_{i-m:n}}{\frac{i+m}{n}-\frac{i-m}{n}}=\frac{X_{i+m:n}-X_{i-m:n}}{2m/n}.  
	\end{align*}
	Analogous to Vasicek (1976), Park (1999), Xiong et al. (2020), and Jose and Sathar (2022), Gupta and Chaudhary (2023), we write The first estimator $\reallywidehat{\Delta_1}$ of $VJ(X)$ as	
	
	\begin{align*}
		\reallywidehat{\Delta_1}=\frac{1}{4n} \sum_{i=1}^{n} \left(\frac{2 m/n}{X_{i+m:n}-X_{i-m:n}}\right)^{2} - \frac{1}{4} \left[ \frac{1}{n} \sum_{i=1}^{n} \left(\frac{2 m/n}{X_{i+m:n}-X_{i-m:n}}\right)\right]^2
	\end{align*}
	
	Here window size $m$  is a positive integer less than $\frac{n}{2}.$  If $i+m > n$ then we consider $X_{i+m:n} = X_{n:n}$ and if $i-m < 1$ then we consider $X_{i-m:n} = X_{1:n}.$ 
	
	The following theorem shows $\reallywidehat{\Delta_1}$ is a consistent estimator of $VJ(X).$ Proof is similar to Vasicek (1976) and hence omitted.
	
	\begin{theorem}
		Let $X_1,X_2, \dots, X_n$ be a random sample from an absolutely continuous cumulative distribution function $F$ and pdf $f.$ Then, $\reallywidehat{\Delta_1}$ converges in probability to $VJ(X)$ as $n \rightarrow \infty, \ \ m \rightarrow \infty, \ \ \frac{m}{n} \rightarrow 0.$
	\end{theorem}
	
	\begin{theorem}\label{thm2}
		Let $X_1, X_2, \dots, X_n$ be a sequence of independent and identically distributed random variables having Varextropy $VJ(X).$ Let $Y_i=aX_i+b, \ a>0, \ b\in R, \ i=1,2,3, \dots, n$ with varextropy $VJ(Y).$  Denote the varextropy estimators based on $X_i$ and  $Y_i$ by $\reallywidehat{\Delta^X_1}$ and $\reallywidehat{\Delta^Y_1}$ respectively. Then, the following properties hold.\\
		\begin{enumerate}
			\item [(i)]$E(\reallywidehat{\Delta^Y_1})=\frac{1}{a^2}E(\reallywidehat{\Delta^X_1})$
			\item [(ii)] $Var(\reallywidehat{\Delta^Y_1})=\frac{1}{a^4} Var(\reallywidehat{\Delta^X_1})$
			\item [(iii)] $MSE(\reallywidehat{\Delta^Y_1})=\frac{1}{a^4} MSE(\reallywidehat{\Delta^X_1})$
		\end{enumerate}
	\end{theorem}
	\noindent \textbf{Proof} We have
	\begin{align*}
		\reallywidehat{\Delta^Y_1}&= \frac{1}{4n} \sum_{i=1}^{n} \left(\frac{2 m/n}{Y_{i+m:n}-Y_{i-m:n}}\right)^{2} - \frac{1}{4} \left[ \frac{1}{n} \sum_{i=1}^{n} \left(\frac{2 m/n}{Y_{i+m:n}-Y_{i-m:n}}\right)\right]^2
		\\
		&=\frac{1}{4n} \sum_{i=1}^{n} \left(\frac{2 m/n}{(aX_{i+m:n}+b)-(aX_{i-m:n}+b)}\right)^{2} - \frac{1}{4} \left[ \frac{1}{n} \sum_{i=1}^{n} \left(\frac{2 m/n}{aX_{i+m:n}+b)-(aX_{i-m:n}+b)}\right)\right]^2 \\
		&= \frac{1}{a^2}\left(\frac{1}{4n} \sum_{i=1}^{n} \left(\frac{2 m/n}{X_{i+m:n}-X_{i-m:n}}\right)^{2} - \frac{1}{4} \left[ \frac{1}{n} \sum_{i=1}^{n} \left(\frac{2 m/n}{X_{i+m:n}-X_{i-m:n}}\right)\right]^2 \right)	\\
		&=\frac{1}{a^2} \reallywidehat{\Delta^X_1}					
	\end{align*}
	Therefore, the proof of this theorem is completed.

	\subsection{The second estimator}
	
	Ebrahimi et al. (1994) modified the estimator proposed by Vasicek (1976) and obtained a new estimator with less bias and less mean squared error than Vasicek's estimator. Analogous to Ebrahimi et al. (1994), we propose an estimator $\reallywidehat{\Delta_2}$ as
	
	\begin{align*}
		\reallywidehat{\Delta_2}&=\frac{1}{4n} \sum_{i=1}^{n} \left(\frac{c_i m/n}{X_{i+m:n}-X_{i-m:n}}\right)^{2} - \frac{1}{4} \left[ \frac{1}{n} \sum_{i=1}^{n} \left(\frac{c_i m/n}{X_{i+m:n}-X_{i-m:n}}\right)\right]^2,
	\end{align*}
	
	where 
	\begin{align*}
		c_i=		
		\begin{cases}
			1+\frac{i-1}{m}, \hspace{10mm} 1\leq i \leq m\\
			2, \hspace{16mm} m+1\leq i \leq n-m\\
			1+\frac{n-i}{m}, \hspace{4mm} n-m+1\leq i \leq n.
		\end{cases}~~~~~~~~~~~~~~~~~
	\end{align*}
	
	The following theorem shows $\reallywidehat{\Delta_2}$ is a consistent estimator of $VJ(X).$ Proof is similar to Ebrahimi et al. (1994) and hence omitted.
	
	\begin{theorem}
		Let $X_1, X_2, \dots, X_n$ be a random sample from an absolutely continuous cumulative distribution function $F$ and pdf $f.$ Then, $\reallywidehat{\Delta_2}$ converges in probability to $VJ(X)$ as $n \rightarrow \infty, \ \ m \rightarrow \infty, \ \ \frac{m}{n} \rightarrow 0.$
	\end{theorem}
	
	\begin{theorem}
		Let $X_1, X_2, \dots, X_n$ be a sequence of independent and identically distributed random variables having Vaextropy $VJ(X).$ Let $Y_i=aX_i+b, \ a>0, \ b\in R, \ i=1,2,3, \dots, n$ with varextropy $VJ(Y).$  Denote the varextropy estimators based on $X_i$ and  $Y_i$ by $\reallywidehat{\Delta^X_2}$ and $\reallywidehat{\Delta^Y_2}$ respectively. Then, the following properties hold.\\
		\begin{enumerate}
			\item [(i)]$E(\reallywidehat{\Delta^Y_2})=\frac{1}{a^2}E(\reallywidehat{\Delta^X_2})$
			\item [(ii)] $Var(\reallywidehat{\Delta^Y_2})=\frac{1}{a^4} Var(\reallywidehat{\Delta^X_2})$
			\item [(iii)] $MSE(\reallywidehat{\Delta^Y_2})=\frac{1}{a^4} MSE(\reallywidehat{\Delta^X_2})$
		\end{enumerate}
	\end{theorem}
	\noindent \textbf{Proof} Proof is similar to the proof of Theorem 
	\ref{thm2}.
	
	\subsection{The third estimator}
	Dmitriew and Tarasenko (1973) proposed to estimate $\int_{-\infty}^{+\infty} f^2(x)dx$ by $\int_{-\infty}^{+\infty} \hat{f}^2(x)dx$. Hence, the varextropy of $X$ can be estimated by $\reallywidehat{\Delta_3}$ as
	\begin{align*}
		\reallywidehat{\Delta_3}=\frac{1}{4}  \int_{S_X} \hat{f}^3(x)dx - \frac{1}{4} \left(\int_{S_X} \hat{f}^2(x)dx\right)^2,
	\end{align*}
	
	where $\hat{f}$ is kernel density function estimator of $f$ which is defined as
	
	\begin{align*}
		\hat{f}(x)=\frac{1}{nh} \sum_{i=1}^{n} k\left(\frac{x-X_i}{h}\right),
	\end{align*}
	where $h$ is the bandwidth and $k$ is a kernel function such that 
	\begin{align*}
		\int_{-\infty}^{\infty} k(x)dx=1.
	\end{align*}
	Usually, $k$  will be a symmetric pdf of some random variable. For example, kernel function $k$ may be considered as pdf of standard normal random variable, that is,
	\[k(x)=\frac{1}{\sqrt{2\pi}} e^{-\frac{1}{2}x^2}, \ \ x \in R.\]

	\subsection{The fourth estimator}
	Since 
	\begin{align*}
		VJ(X)=\frac{1}{4}E(f^2(X))-\frac{1}{4} \left[E(f(X))\right]^2.
	\end{align*}
	One can obtain the simple varextropy estimator as
	
	\begin{align*}
		\reallywidehat{\Delta_4}=\frac{1}{4n} \sum_{i=1}^{n}f^2(X_i)-\frac{1}{4} \left[\frac{1}{n} \sum_{i=1}^{n}f(X_i)\right]^2.
	\end{align*}
	The varextropy $VJ(X)$ of an unknown continuous pdf $f$ can be estimated as
	\begin{align*}
		\reallywidehat{\Delta_4}=\frac{1}{4n} \sum_{i=1}^{n}(\hat{f}^2(X_i))-\frac{1}{4} \left[\frac{1}{n} \sum_{i=1}^{n}\hat{f}(X_i)\right]^2, 
	\end{align*}
	where 
	\begin{align*}
		\hat{f}(X_i)=\frac{1}{nh} \sum_{j=1}^{n} k\left(\frac{X_i-X_j}{h}\right), \  k(x)=\frac{1}{\sqrt{2\pi}} e^{-x^2/2},\ x \in R \ \text{and} \ h=1.06sn^{-1/5}.
	\end{align*}
	The bandwidth $h$ is taken as a normal smoothing formula in which $s$ is the sample standard deviation.

	\begin{theorem}
		Let $X_1, X_2, \dots, X_n$ be a random sample from an absolutely continuous cumulative distribution function $F$ and pdf $f$ with finite variance. Then, $\reallywidehat{\Delta_4}$ converges in probability to $VJ(X)$ as $n \rightarrow \infty.$
	\end{theorem}
	\noindent \textbf{Proof}  Proof follows from the kernel density estimator consistency property and sample mean's consistency property.

	\begin{theorem}
		Let $X_1, X_2, \dots, X_n$ be a sequence of independent and identically distributed random variables having Vaextropy $VJ(X).$ Let $Y_i=aX_i+b, \ a>0, \ b\in R, \ i=1,2,3, \dots, n$ with varextropy $VJ(Y).$  Denote the varextropy estimators based on $X_i$ and  $Y_i$ by $\reallywidehat{\Delta^X_4}$ and $\reallywidehat{\Delta^Y_4}$ respectively. Then, the following properties hold.\\
		\begin{enumerate}
			\item [(i)]$E(\reallywidehat{\Delta^Y_4})=\frac{1}{a^4}E(\reallywidehat{\Delta^X_4})$
			\item [(ii)] $Var(\reallywidehat{\Delta^Y_4})=\frac{1}{a^4} Var(\reallywidehat{\Delta^X_4})$
			\item [(iii)] $MSE(\reallywidehat{\Delta^Y_4})=\frac{1}{a^4} MSE(\reallywidehat{\Delta^X_4})$
		\end{enumerate}
	\end{theorem}
	\noindent \textbf{Proof}  Since $\hat{f}(Y_i)=\frac{1}{a} \hat{f}(X_i),$ and $h_Y=ah_X.$ Therefore, $\reallywidehat{\Delta^Y_4}=\frac{1}{a^2} \reallywidehat{\Delta^X_4}$ and the proof is complete.

	\subsection{The fifth estimator}
	Correa (1995) proposed an estimator based on local linear regression. Analogous to  Correa's estimator, we propose a varextropy estimator based on the following local linear model. 
	Using sample infomation  $((X_{1:n}, F_n(X_{1:n})), (X_{2:n}, F_n(X_{2:n})), \dots, (X_{n:n}, F_n(X_{n:n}))),$ we write 	
	\begin{align*}
		T=\frac{1}{4n}  \sum_{i=1}^{n} \left( \frac{(i+m)/n-(i-m)/n}{X_{i+m:n}-X_{i-m:n}} \right)^3 - \frac{1}{4} \left[\frac{1}{n} \sum_{i=1}^{n} \left( \frac{(i+m)/n-(i-m)/n}{X_{i+m:n}-X_{i-m:n}} \right) \right]^2.
	\end{align*}
	Note that $\frac{(i+m)/n-(i-m)/n}{X_{i+m:n}-X_{i-m:n}}$ is the slope of line joining points $(X_{i+m}, F_n(X_{i+m}))$ and $(X_{i-m}, F_n(X_{i-m})).$ In order to calculate the density of $F(x)$ in the range $(X_{(i-m)}, X_{(i+m))},$ we use a local linear model based on $2m+1$ points,
	
	\[F(x_j)=\alpha + \beta x_{(j)} + \epsilon, \ \ j=m-i, \dots, m+i.\]
	
	The estimator of varextropy is obtained using the least square method as
	
	\begin{align*}
		\Delta_5=\frac{1}{4n}  \sum_{i=1}^{n} \left( \frac{\sum_{j=i-m}^{i+m} (X_{j:n}-\bar{X}_{i:n})(j-i)}{n\sum_{j=i-m}^{i+m} (X_{j:n}-\bar{X}_{i:n})^2}  \right)^3 - \frac{1}{4} \left[\frac{1}{n} \sum_{i=1}^{n} \left( \frac{\sum_{j=i-m}^{i+m} (X_{j:n}-\bar{X}_{i:n})(j-i)}{n\sum_{j=i-m}^{i+m} (X_{j:n}-\bar{X}_{i:n})^2}  \right) \right]^2,
	\end{align*}
	where $\bar{X}_{i:n}=\frac{1}{2m+1}\sum_{j=i-m}^{i+m} X_{j:n}.$
	
	The following theorem shows $\reallywidehat{\Delta_5}$ is a consistent estimator of $VJ(X).$ Proof is similar to Correa (1995) and hence omitted.
	
	\begin{theorem}
		Let $X_1, X_2, \dots, X_n$ be a random sample from an absolutely continuous cumulative distribution function $F$ and pdf $f$ with finite variance. Then, $\reallywidehat{\Delta_5}$ converges in probability to $VJ(X)$ as $n \rightarrow \infty.$
	\end{theorem}

	\begin{theorem}
		Let $X_1, X_2, \dots, X_n$ be a sequence of independent and identically distributed random variables having Vaextropy $VJ(X).$ Let $Y_i=aX_i+b, \ a>0, \ b\in R, \ i=1,2,3, \dots, n$ with varextropy $VJ(Y).$  Denote the varextropy estimators based on $X_i$ and  $Y_i$ by $\reallywidehat{\Delta^X_5}$ and $\reallywidehat{\Delta^Y_5}$ respectively. Then, the following properties hold.\\
		\begin{enumerate}
			\item [(i)]$E(\reallywidehat{\Delta^Y_5})=\frac{1}{a^4}E(\reallywidehat{\Delta^X_5})$
			\item [(ii)] $Var(\reallywidehat{\Delta^Y_5})=\frac{1}{a^4} Var(\reallywidehat{\Delta^X_5})$
			\item [(iii)] $MSE(\reallywidehat{\Delta^Y_5})=\frac{1}{a^4} MSE(\reallywidehat{\Delta^X_5})$
		\end{enumerate}
	\end{theorem}
	\noindent \textbf{Proof}  Proof is similar to the proof of Theorem 
	\ref{thm2}.

	\subsection{The sixth estimator}
	
	Noughabi and Noughabi (2023) introduced some estimators of the varentropy of a continuous random variable. Analogous to Noughabi and Noughabi (2023), we introduce here an estimator of the varextropy of a continuous random variable.
	
	By using \[\frac{d}{dp}(F^{-1}(p)) \approxeq \frac{X_{i+m}-X_{i-m}}{F(X_{i+m})-F(X_{i-m})}, \]
	
	for $\frac{i-1}{n}<p\leq \frac{i}{n}, \ i=m+1, m+2, \dots, n-m, $ where $m$ is positive integer smaller than $n/2.$ The vaextropy of X can be written as 
	
	\begin{align*}
		VJ(X)&=\frac{1}{4} \int_{0}^{1} \left(\frac{d}{dp}(F^{-1}(p))\right)^{-2} dp- \frac{1}{4} \left[\int_{0}^{1} \left( \frac{d}{dp}(F^{-1}(p))\right)^{-1}dp\right]^2 \\
		&\approxeq \frac{1}{4} \sum_{i=1}^{n} \left( \frac{X_{i+m:n}-X_{i-m:n}}{F(X_{i+m:n})-F(X_{i-m:n})}\right)^{-2} (p_{i+1}-p_{i}) \\
		& \ \ -\frac{1}{4} \left[\sum_{i=1}^{n} \left( \frac{X_{i+m:n}-X_{i-m:n}}{F(X_{i+m:n})-F(X_{i-m:n})} \right)^{-1} (p_{i+1}-p_{i}) \right]^2
	\end{align*}
	Now, we propose approximate equality as follows. 
	\[F(X_{i+m:n})- F(X_{i-m:n}) \approxeq \frac{1}{2} (f(X_{i+m:n})- f(X_{i-m:n}))(X_{i+m:n}- X_{i-m:n}) \]
	
	Therefore, we have
	\begin{align*}
		VJ(X)&\approxeq \frac{1}{4n} \sum_{i=1}^{n} \left( \frac{X_{i+m:n}-X_{i-m:n}}{\frac{1}{2} (f(X_{i+m:n})- f(X_{i-m:n}))(X_{i+m:n}- X_{i-m:n})}\right)^{-2}\\
		& \ \ \ \ - \frac{1}{4} \left[ \frac{1}{n} \sum_{i=1}^{n} \left( \frac{X_{i+m:n}-X_{i-m:n}}{\frac{1}{2} (f(X_{i+m:n})- f(X_{i-m:n}))(X_{i+m:n}- X_{i-m:n})} \right)^{-1} \right]^2\\
		&= \frac{1}{4n} \sum_{i=1}^{n} \left( \frac{2}{ (f(X_{i+m:n})- f(X_{i-m:n}))}\right)^{-2}\\
		& \ \ \ \ - \frac{1}{4} \left[ \frac{1}{n} \sum_{i=1}^{n} \left( \frac{2}{ (f(X_{i+m:n})- f(X_{i-m:n}))} \right)^{-1} \right]^2 \\    
		&= \frac{1}{4n} \sum_{i=1}^{n} \left( \frac{(f(X_{i+m:n})- f(X_{i-m:n}))}{2 }\right)^{2}\\
		& \ \ \ \ - \frac{1}{4} \left[ \frac{1}{n} \sum_{i=1}^{n} \left( \frac{(f(X_{i+m:n})- f(X_{i-m:n}))}{2} \right) \right]^2    
	\end{align*}
	The varextropy $VJ(X)$ of an unknown continuous pdf f is therefore proposed to be estimated by
	
	\begin{align*}
		\Delta_6&= \frac{1}{4n} \sum_{i=1}^{n} \left( \frac{(\hat{f}(X_{i+m:n})- \hat{f}(X_{i-m:n}))}{2 }\right)^{2}\\
		& \ \ \ \ - \frac{1}{4} \left[ \frac{1}{n} \sum_{i=1}^{n} \left( \frac{(\hat{f}(X_{i+m:n})- \hat{f}(X_{i-m:n}))}{2} \right) \right]^2,   
	\end{align*}
	where 
	\[\hat{f}(X_{i})= \frac{1}{nh} \sum_{j=1}^{n} k\left(\frac{X_i-X_j}{h}\right), \]
	and the bandwidth h is chosen to be the normal optimal smoothing formula $h=1.06 sn^{-1/5}$in which s is the sample standard deviation and the kernel function is chosen to be the standard normal density function. We will consider $X_{(i)}=X_{(1)}$ if $i<1$ and $X_{(i)}=X_{(n)}$ if $i>n.$
	
	The following theorem shows $\reallywidehat{\Delta_5}$ is a consistent estimator of $VJ(X).$ Proof is similar to the proof of Theorem 2.8 of Noughabi and Noughabi (2023) and hence omitted.

	\begin{theorem}
		Let $X_1, X_2, \dots, X_n$ be a random sample from an absolutely continuous cumulative distribution function $F$ and pdf $f$ with finite variance. Then, $\reallywidehat{\Delta_6}$ converges in probability to $VJ(X)$ as $n \rightarrow \infty.$
	\end{theorem}

	\begin{theorem}
		Let $X_1, X_2, \dots, X_n$ be a sequence of independent and identically distributed random variables having Vaextropy $VJ(X).$ Let $Y_i=aX_i+b, \ a>0, \ b\in R, \ i=1,2,3, \dots, n$ with varextropy $VJ(Y).$  Denote the varextropy estimators based on $X_i$ and  $Y_i$ by $\reallywidehat{\Delta^X_6}$ and $\reallywidehat{\Delta^Y_6}$ respectively. Then, the following properties hold.\\
		\begin{enumerate}
			\item [(i)]$E(\reallywidehat{\Delta^Y_6})=\frac{1}{a^4}E(\reallywidehat{\Delta^X_6})$
			\item [(ii)] $Var(\reallywidehat{\Delta^Y_6})=\frac{1}{a^4} Var(\reallywidehat{\Delta^X_6})$
			\item [(iii)] $MSE(\reallywidehat{\Delta^Y_6})=\frac{1}{a^4} MSE(\reallywidehat{\Delta^X_6})$
		\end{enumerate}
	\end{theorem}
	\noindent \textbf{Proof}  Proof is similar to the proof of Theorem 
	\ref{thm2}.

	\section{A Characterization of Uniform distribution}
	In many practical problems, the goodness-of-fit test may be reduced to the problem of testing uniformity.  Since varextropy of $X$ is the variance of $-\frac{1}{2}f(X),$, therefore, varextropy is a non-negative for any random variable $X.$ Among all distributions with support on [0,1], the uniform distribution has the maximum extropy. An important property of uniform distribution is that it obtains the minimum varextropy among all distributions having support on [0,1] (see Qiu and Jia (2017)). 
	\begin{theorem}\label{uniformunique}
		Let X be a continuous random variable having support on [0,1].	Then $VJ(X)= 0$ if and only if X has a uniform distribution on the interval [0,1].
	\end{theorem}
	\noindent \textbf{Proof} Let random variable X have a uniform distribution on the interval [0,1], then $f(x)=1,\ 0 \leq x \leq 1$ and 
	\begin{align*} 
		VJ(X)=\frac{1}{4}  \int_{0}^{1} f^3(x)dx - \frac{1}{4} \left[ \int_{0}^{1} f^2(x)dx \right]^2	=0
	\end{align*}
	
	Conversely, $VJ(X)= 0$ implies $Var(f(X))=0,$ that is, $f(x)=c$. Since 
	$\int_{0}^{1} f(x) dx=1,$ therefore $f(x)=1, \ 0\leq x \leq 1.$ Hence proof is complete.

	\section{Test of symmetry using extropy }
	Gupta and Chaudhary (2023) proved the random variable $X$ has a symmetric distribution if and only if for a fixed $k\geq 1$,  $\bar{\xi} J(L_{n,k})=\xi J(U_{n,k})$ for all $n\geq 1$. that will be used to derive test statistics for testing the symmetry of a continuous symmetric distribution. This motivates us that ${\Delta}_{n,k} =\xi J(U_{n,k})-\bar{\xi}  J(L_{n,k})$ can be used to test whether $X$ is symmetric distribution. We reject the null hypothesis of symmetry because small or large values of ${\Delta}_{n,k}$ can be seen as symptoms of non-symmetry. For this reason, we propose a test for symmetry based on the sample estimator of ${\Delta}_{n,k}.$ We see that  ${\Delta}_{n,k} =0$ if and only if $X$ is a symmetric distribution. Therefore, if an independent and identically distributed (iid) sample of size $N$ is available, its empirical counterpart ${\reallywidehat{\Delta}}_{n,k}$ will help in determining whether or not the sample comes from a symmetric distribution.
	
	\noindent Cumulative past extropy  of $n$th lower $k$-record value is 
	\begin{align*}
		\bar{\xi} J (L_{n,k}) &=- \frac{1}{2} \int_{S_X} F_{L_{n,k}}^2 (x) dx = - \frac{1}{2} \int_{0}^{1} u^{2k} \left( \sum_{j=0}^{n-1} \frac{(-k\log (u))^j }{j!}  \right)^2 \left(\frac{d}{du} F^{-1}(u) \right)du.
	\end{align*}
	
	\noindent Cumulative residual extropy  of $n$th upper $k$-record value is 
	\begin{align}
		\xi J (U_{n,k}) &=- \frac{1}{2} \int_{S_X} \bar{F}_{U_{n,k}}^2 (x) dx = - \frac{1}{2} \int_{0}^{1} (1-u)^{2k} \left( \sum_{j=0}^{n-1} \frac{(-k\log (1-u))^j  }{j!}  \right)^2 \frac{d}{du} (F^{-1}(u))du. \nonumber
	\end{align}
	
	\noindent We propose test statistics as:
	\begin{align*}
		{{\Delta}}_{n,k}&=\xi J (U_{n,k})-\bar{\xi} J (L_{n,k}) \\
		&=  - \frac{1}{2} \int_{0}^{1}  \left[(1-u)^{2k} \left( \sum_{j=0}^{n-1} \frac{(-k\log (1-u))^j }{j!}  \right)^2 \right.  \nonumber \\ & \ \ \ \ \ \ \ \ \  \ \ \ \ \ \ \ \   \ \ \ \ \ \ \ \  \ \ \left. - u^{2k} \left( \sum_{j=0}^{n-1} \frac{(-k\log (u))^j u^k}{j!}  \right)^2   \right] \left(\frac{dF^{-1} (u)}{du}\right) du .
	\end{align*}
	
	\noindent An estimator of test statistics ${{\Delta}}_{n,k}$ is 
	\begin{align}
		&{\reallywidehat{\Delta}}_{n,k} =\reallywidehat{\xi J (U_{n,k})}-\reallywidehat{\bar{\xi} J (L_{n,k})} \nonumber\\
		&= - \frac{1}{2N} \sum_{i=1}^{N} \left[\left(1-\frac{i}{N+1}\right)^{2k} \left( \sum_{j=0}^{n-1} \frac{(-k\log (1-\frac{i}{N+1}))^j  }{j!}  \right)^2 \right. \nonumber\\ &\  \left.-\left(\frac{i}{N+1}\right)^{2k} \left( \sum_{j=0}^{n-1} \frac{(-k\log (\frac{i}{N+1}))^j }{j!}  \right)^2   \right] \frac{(X_{i+m:N}-X_{i-m:N})}{2m/N} .\nonumber 
	\end{align}		
	
	We choose $n=2$ and $k=2$ for the simplicity of calculation, so we consider test statistics as $\reallywidehat{\Delta}_{2,2} $. 
	\begin{align}
		{\reallywidehat{\Delta}}_{2,2} &= - \frac{1}{2N} \sum_{i=1}^{N} \left[\left(1-\frac{i}{N+1}\right)^{4} \left(1-2\log (1-\frac{i}{N+1})  \right)^2\right.\nonumber\\ & \ \ \ \ \  \left. -\left(\frac{i}{N+1}\right)^{4} \left( 1-2\log (\frac{i}{N+1})\right)^2   \right] \frac{(X_{i+m:N}-X_{i-m:N})}{2m/N}. \nonumber
	\end{align}		
	
	The procedure remains the same when one chooses any other $n$ and $k$. Following theorem says ${\reallywidehat{\Delta}}_{2,2}$ is consistent estimator of ${\Delta_{2,2}}.$ 
	\begin{theorem}
		Assume that $X_1,\ X_2,\,..., X_N$ is a random sample of size $N$ taken from a population with pdf $f$ and cdf $F$. Also, let the variance of the random variable be finite. Then ${\reallywidehat{\Delta}}_{2,2}$ converges in probability to $\Delta_{2,2}$ as $N \longrightarrow \infty, \ m\longrightarrow \infty \ \text{and} \ \frac{m}{N} \longrightarrow 0.$
	\end{theorem}
	\textbf{Proof : } Following lines of the proof of Theorem 1 of Vasicek (1976). We have, $\reallywidehat{\xi J (U_{2,2})} \overset{P}{\to} \xi J (U_{2,2})$ and  $\reallywidehat{\bar{\xi}J (L_{2,2})} \overset{P}{\to} \bar{\xi}J (L_{2,2})$. Therefore we get,  ${\reallywidehat{\Delta}}_{2,2} \overset{P}{\to}  {\Delta}_{2,2}$. That is, ${\reallywidehat{\Delta}}_{2,2}$ is consistent estimator of  ${\Delta}_{2,2}$.\\
	
	\begin{theorem}
		Let $X_1, X_2, ... , X_N$ be a sequence of iid random variables and let $Y_i=aX_i+b, \ a>0,\ b\in \mathbb{R},\ i=1,2,... ,N.$ Denote the estimator for $\Delta_{2,2}$ based on $X_i$ and $Y_i$ by $\reallywidehat{\Delta_{2,2}^X}$ and $\reallywidehat{\Delta_{2,2}^Y}$, respectively. Then 
		\begin{enumerate}[(i)] 
			\item  E($\reallywidehat{\Delta_{2,2}^Y})=a E(\reallywidehat{\Delta_{2,2}^X})$
			\item  Var($\reallywidehat{\Delta_{2,2}^Y})=a^2 Var(\reallywidehat{\Delta_{2,2}^X})$
			\item  MSE($\reallywidehat{\Delta_{2,2}^Y})=a^2 MSE(\reallywidehat{\Delta_{2,2}^X})$\\
		\end{enumerate}
		where $E(X),\ Var(X)$ and $MSE(X)$ represent expectation, variance and mean square error of random variable $X$, respectively.
	\end{theorem}
	\textbf{Proof}
	\begin{align}
		\newline   \reallywidehat{{\Delta}_{2,2}^Y }
		= & - \frac{1}{2N} \sum_{i=1}^{N} \left[\left(1-\frac{i}{N+1}\right)^{4}  \left(1-2\log ( 1 - \frac{i}{N+1} ) \right)^2 \right.\nonumber \\ & \left.-\left(\frac{i}{N+1}\right)^{4} \left( 1-2\log (\frac{i}{N+1})\right)^2   \right] \frac{(Y_{i+m:N}-Y_{i-m:N})}{2m/N} \nonumber\\
		= & - \frac{1}{2N} \sum_{i=1}^{N} \left[\left(1-\frac{i}{N+1}\right)^{4} \left(1-2\log (1-\frac{i}{N+1})  \right)^2  \right.\nonumber \\ & \left. -\left(\frac{i}{N+1}\right)^{4} \left( 1-2\log (\frac{i}{N+1})\right)^2   \right] \frac{(aX_{i+m:N}-aX_{i-m:N})}{2m/N} \nonumber\\
		= & a\reallywidehat{{\Delta}_{2,2}^X }.\nonumber
	\end{align}
	Thus, proof is completed because of $\reallywidehat{{\Delta}_{2,2}^Y}=a \reallywidehat{{\Delta}_{2,2}^X}$ and  properties of mean, variance and MSE of $X$.
	\subsection{Critical values}
	For the computation of the critical values, it is crucial to assess the asymptotic distribution of the newly suggested test statistic. Unfortunately, it is very difficult to deduce the distribution of ${\reallywidehat{\Delta}}_{2,2}$ as $N\rightarrow$ $\infty$, because there is a window size $m$ that depends on $N$. 
	
	Now, 10,000 samples of size $N$ taken from the standard normal distribution, are used to calculate the critical values of the ${\reallywidehat{\Delta}}_{2,2}$ function. Then, the values of ${\reallywidehat{\Delta}}_{2,2}$ are obtained based on these 10,000 samples of size $N$ generated from the standard normal distribution. From the 10,000 values of  ${\reallywidehat{\Delta}}_{2,2}$, $(1-\frac{\alpha}{2})$th quantile represents the critical value corresponding to sample size $N$ of the test statistic at significance level $\alpha$. If the critical values are denoted as ${\reallywidehat{\Delta}}_{2,2}(1-\frac{\alpha}{2})$, then the null hypothesis is rejected with size $\alpha$ whenever $|{\reallywidehat{\Delta}}_{2,2}|>{\reallywidehat{\Delta}}_{2,2}(1-\frac{\alpha}{2}).$

	The exact critical values of ${\reallywidehat{\Delta}}_{2,2}$ based on 10,000 samples of different sizes generated from the standard normal distribution at significance level  $\alpha=0.05$ is given in Table 1. The critical values are obtained for sample sizes $N=5,10,20,30,40,50,100$ with window sizes $m$ ranging from 2 to 40. The next section deals with the simulation study through which the power of the test statistic is evaluated. 
	
	\begin{center}
		{\bf Table 1}.  \small{{Critical values of $|{\reallywidehat{\Delta}}_{2,2}|$ statistics at significance level $\alpha$= 0.05}} 
		
		\resizebox{!}{!}{
			\begin{tabular} {p{1.0cm} p{1.0cm} p{1.0cm} p{1.0cm} p{1.0cm} p{1.0cm} p{1.0cm} p{1.0cm} p{1.0cm} p{1.0cm} p{1.0cm} } 
				\hline
				$m\backslash N$ & 5 & 10 & 20  & 30  & 40 & 50 & 100 \\
				\hline \\
				2 & 0.3637 & 0.6093  & 0.6673  & 0.6703 & 0.6658 & 0.6474 & 0.5969 \\
				3 & \   & 0.4787  & 0.5833  & 0.5936 & 0.6011 & 0.5857 & 0.5611 \\ 
				4 & \   & 0.3641  & 0.5333  & 0.5539 & 0.5553 & 0.5387 & 0.5284 \\
				5 & \   & \  & 0.4776 & 0.5216 & 0.5287 & 0.5305 & 0.5054 \\
				6 & \   & \  & 0.4362  & 0.4872 & 0.5074 & 0.4979 & 0.4794 \\
				7 & \   & \  & 0.3848  & 0.4536 & 0.4785 & 0.4718 & 0.4750 \\
				8 & \   & \  & 0.3460  & 0.4207 & 0.4543 & 0.4647 & 0.4573 \\
				9  & \   & \  & 0.2951  & 0.4044 & 0.4271 & 0.4491 & 0.4488 \\
				10 & \   & \  & \  & 0.3642 &  0.4188 & 0.4235 & 0.4405 \\
				11 & \   & \  & \  & 0.3440 & 0.3953 & 0.4259 & 0.4242 \\
				12 & \   & \  & \  & 0.3254 & 0.3775 & 0.4053 & 0.4249 \\
				13 & \   & \  & \  & 0.2948 & 0.3561 & 0.3947 & 0.4251 \\
				14 & \   & \  & \  & 0.2751 & 0.3515 & 0.3821 & 0.4203 \\
				15 & \   & \  & \  & \ & 0.3239 & 0.3626 & 0.4157 \\
				16 & \   & \  & \  & \ & 0.3057 & 0.3521 & 0.3967  \\
				17 & \   & \  & \  & \ & 0.2868 & 0.3433 & 0.3901 \\
				18 & \   & \  & \  & \ & 0.2769 & 0.3356 & 0.3820 \\
				19 & \   & \  & \  & \ & 0.2583 & 0.3182 & 0.3883 \\
				20 & \   & \  & \  & \ & \ & 0.3011 & 0.3779 \\
				21 & \   & \  & \  & \ & \ & 0.2909 & 0.3646 \\
				22 & \   & \  & \  & \ & \ & 0.2777 & 0.3705 \\
				23 & \   & \  & \  & \ & \ & 0.2697 & 0.3607 \\
				24 & \   & \  & \  & \ & \ & 0.2533 & 0.3575 \\
				25 & \   & \  & \  & \ & \ & \ & 0.3592 \\
				26 & \   & \  & \  & \ & \ & \ & 0.3446 \\
				27 & \   & \  & \  & \ & \ & \ & 0.3358 \\
				28 & \   & \  & \  & \ & \ & \ & 0.3427 \\
				29 & \   & \  & \  & \ & \ & \ & 0.3355 \\
				30 & \   & \  & \  & \ & \ & \ & 0.3258 \\
				40 & \   & \  & \  & \ & \ & \ & 0.2797 \\

				\hline
		\end{tabular}}
	\end{center}\label{table1}

	\subsection{Power and size of test}
	Unfortunately, it is very complicated to derive the exact distribution of ${\reallywidehat{\Delta}}_{2,2}$ because it depends on window size $m$ which is dependent on sample size $N$. Tables 1, 2 and 3 show the exact critical values of the test statistic $|{\reallywidehat{\Delta}}_{2,2}|$ for various sample sizes by Monte Carlo simulation with 10,000 repetitions, for significance levels $\alpha = 0.05.$ A similar procedure has been used in Xiong et al. (2020). To determine whether the test statistic's absolute value is greater than the critical value, we generated a sample of size $N$ from the standard normal distribution. We then repeated this process 10,000 times. The power of the test is measured by the percentage of rejection. Tables 4, 5 and 6 give the power of the test when an alternative distribution is taken as a chi-square distribution with a degree of freedom 1 for $\alpha= 0.10$,  $\alpha = 0.05$ and $\alpha = 0.01$, respectively. Since $\chi^2(1)$ is asymmetric distribution. Our test significantly verifies this fact. We observe from Tables 4, 5 and 6 that when we increase sample size $N$, power increases. When the sample size is 100, power is 1.000 for any value of $m$. That means our test performs well for a large sample size. 	Table 7 gives information about power against different alternative distribution $\chi^2{(1)}$, $\chi^2{(2)}$, $\chi^2{(3)}$ and $N(0,1)$. Since N(0,1) is a symmetric distribution,  therefore power against null distribution is the size of the test. 	The powers of ${\reallywidehat{\Delta}}_{2,2}$  should be near a significant level if the alternative distribution is symmetric. Tables 7 and 8 use the standard normal distribution as an example, and we can see that all of the powers are roughly equal to 0.05 as they should be in this situation. As a result, ${\reallywidehat{\Delta}}_{2,2}$  can maintain its nominal level. A smaller value of $m$  works better for any sample size. Table \ref{table9} lists a proposed value of window size $m$ for different sample sizes based on the simulation study.
	
	\begin{center}
		\noindent {\bf Table 2}.  \small{{Powers of  ${\reallywidehat{\Delta}}_{2,2}$ statistics against alternative $\chi^2_{(1)}$ at significance level $\alpha$= 0.05}} \\ 
		\resizebox{!}{!}{
			\begin{tabular}{ p{1.0cm} p{1.0cm} p{1.0cm} p{1.0cm} p{1.0cm} p{1.0cm} p{1.0cm} p{1.0cm} } 
				\hline
				$m\backslash N$ & 5 & 10 & 20  & 30  & 40 & 50 & 100 \\
				\hline \\
				2  & 0.2649   & 0.5644  & 0.8759  & 0.9685 & 0.9924 & 0.9984 & 1.0000 \\
				3 & \   & 0.5624  & 0.8859  & 0.9703 & 0.9922 & 0.9991 & 1.0000 \\ 
				4 & \   & 0.5534  & 0.8769  & 0.9749 & 0.9957 & 0.9997 & 1.0000 \\
				5 & \   & \  & 0.8794  & 0.9694 & 0.9952 & 0.9990 & 1.0000 \\
				6 & \   & \  & 0.8691  & 0.9784 & 0.9941 & 0.9991 & 1.0000 \\
				7 & \   & \  & 0.8756  & 0.9756 & 0.9945 & 0.9995 & 1.0000 \\
				8 & \   & \  & 0.8536  & 0.9759 & 0.9952 & 0.9989 & 1.0000 \\
				9  & \   & \  & 0.8556  & 0.9680 & 0.9949 & 0.9993 & 1.0000 \\
				10 & \   & \  & \  & 0.9679 & 0.9936 & 0.9991 & 1.0000 \\
				11 & \   & \  & \  & 0.9708 & 0.9937 & 0.9986 & 1.0000 \\
				12 & \   & \  & \  & 0.9609 & 0.9928 & 0.9991 & 1.0000 \\
				13 & \   & \  & \  & 0.9589 & 0.9945 & 0.9991 & 1.0000 \\
				14 & \   & \  & \  & 0.9586 & 0.9913 & 0.9989 & 1.0000 \\
				15 & \   & \  & \  & \ & 0.9924 & 0.9991 & 1.0000 \\
				16 & \   & \  & \  & \ & 0.9902 & 0.9988 & 1.0000 \\
				17 & \   & \  & \  & \ & 0.9902 & 0.9978 & 1.0000 \\
				18 & \   & \  & \  & \ & 0.9903 & 0.9972 & 1.0000 \\
				19 & \   & \  & \  & \ & 0.9870 & 0.9984 & 1.0000 \\
				20 & \   & \  & \  & \ & \ & 0.9987 & 1.0000 \\
				21 & \   & \  & \  & \ & \ & 0.9977 & 1.0000 \\
				22 & \   & \  & \  & \ & \ & 0.9971 & 1.0000 \\
				23 & \   & \  & \  & \ & \ & 0.9964 & 1.0000 \\
				24 & \   & \  & \  & \ & \ & 0.9958 & 1.0000 \\
				25 & \   & \  & \  & \ & \ & \ & 1.0000 \\
				26 & \   & \  & \  & \ & \ & \ & 1.0000 \\
				27 & \   & \  & \  & \ & \ & \ & 1.0000 \\
				28 & \   & \  & \  & \ & \ & \ & 1.0000 \\
				29 & \   & \  & \  & \ & \ & \ & 1.0000 \\
				30 & \   & \  & \  & \ & \ & \ & 1.0000 \\
				40 & \   & \  & \  & \ & \ & \ & 1.0000 \\
				
				\hline
		\end{tabular}}
	\end{center}\label{table2}
	
	\begin{center}
		\noindent{\bf Table 7}.  \small{{Power of ${\reallywidehat{\Delta}}_{2,2}$ for $N=20,50,100$ and significance level $\alpha= 0.05$ for alternatives  $\chi^2_{(1)}$,  $\chi^2_{(2)}$, $\chi^2_{(3)}$  and $N(0,1)$ }} \\
		
		\resizebox{!}{!}{
			\begin{tabular}{ p{1.0cm} p{1.0cm} p{1.0cm} p{1.0cm} p{1.0cm}  p{1.0cm}} 
				\hline
				N & m  & $\chi^2_{(1)}$ & $\chi^2_{(2)}$ & $\chi^2_{(3)}$ & $N(0,1)$ \\
				\hline 
				
				\ & 2   & 0.8759  & 0.8861  & 0.8627 & 0.0518   \\
				\ & 3   & 0.8859  & 0.8976 & 0.8728 & 0.0481  \\
				\ & 4   & 0.8769  & 0.8956 & 0.8762 & 0.0519  \\
				\ & 5   & 0.8794  & 0.8943 & 0.8715 & 0.0478  \\
				20 & 6   & 0.8691  & 0.8941 & 0.8732 & 0.0502 \\
				\ & 7  & 0.8756  & 0.8889 & 0.8716 & 0.0489  \\
				\ & 8   & 0.8536  & 0.8814 & 0.8643 & 0.0491   \\
				\ & 9   & 0.8556  & 0.8757 &0.8555 & 0.0487   \\

				\hline
				
				\ & 2   & 0.9984  & 0.9981 & 0.9955 & 0.0507   \\
				\ & 4   & 0.9997  & 0.9991 & 0.9969 & 0.0516  \\
				\ & 7   & 0.9995  & 0.9987 & 0.9964  & 0.0521   \\
				\ & 9   & 0.9993  & 0.9986 & 0.9960 & 0.0494  \\
				50 & 15   &  0.9991  & 0.9985 & 0.9961 & 0.0466  \\
				\ & 17  & 0.9978  & 0.9986 & 0.9955 & 0.0502  \\
				\ & 20  &  0.9987  & 0.9979 & 0.9944 & 0.0488  \\
				\ & 22  & 0.9971  & 0.9976 & 0.9946 & 0.0497  \\

				\hline
				
				\ & 2   & 1.0000  & 1.0000 & 1.0000 & 0.0496   \\
				\ & 4   & 1.0000  & 1.0000 & 1.0000 & 0.0512  \\
				\ & 5   & 1.0000  & 1.0000 & 1.0000 & 0.0494  \\
				\ & 7   & 1.0000  & 1.0000 & 1.0000 & 0.0515  \\
				100 & 10   & 1.0000  & 1.0000 & 1.0000 & 0.0502   \\
				\ & 15 & 1.0000  & 1.0000 &1.0000 & 0.0473   \\
				\ & 20   & 1.0000  & 1.0000 & 1.0000 & 0.0524   \\
				\ & 30  & 1.0000  & 1.0000 & 1.0000  &  0.0530  \\
				\ & 40   & 1.0000  & 1.0000 & 1.0000 & 0.0520   \\
				
				\hline
		\end{tabular}}
	\end{center}\label{table7}

	\begin{center}
		\noindent{\bf Table 8}.  \small{{Size of ${\reallywidehat{\Delta}}_{2,2}$ for $N=20,50,100$ and significance level $\alpha= 0.05$ for  standard normal distribution. }} \\		
		\vspace{0.2cm}
		\resizebox{!}{!}{
			\begin{tabular}{ p{0.5cm} p{0.5cm} p{1.5cm} p{0.5cm} p{0.5cm} p{1.5cm} p{0.5cm} p{0.5cm} p{1.5cm}} 
				\hline
				N   &m       & $N(0,1)$ 	  &N & m   & $N(0,1)$   &N & m   & $N(0,1)$ \\
				\hline 
				\ & 2       & 0.0518          &\ & 2    & 0.0507  &\ & 2  & 0.0496  \\
				\ & 3       & 0.0481          &\ & 3    & 0.0467  &\ & 3  & 0.0494  \\
				\ & 4       & 0.0519          &\ & 5    & 0.0491  &\ & 5  & 0.0494  \\
				\ & 5       & 0.0478          &\ & 8    & 0.0514  &\ & 8  & 0.0484   \\
				20 & 6      & 0.0502          &50 & 10  & 0.0513  &100 & 10   & 0.0502  \\
				\ & 7       & 0.0489          &\ & 15   & 0.0466  &\ & 15     & 0.0473  \\
				\ & 8       & 0.0491          &\ & 20   & 0.0488  &\ & 20     & 0.0524  \\	
				\ & 9       & 0.0487          &\ & 24   & 0.0503  &\ & 30     & 0.0530  \\	
				\ &\         &\		          &\ &\      &\       &\ & 49     & 0.0515  \\	    
				\hline

		\end{tabular}}
	\end{center}\label{table8}

	\begin{center}
		\noindent{\bf Table 9}.\small{{ Proposed value of window size $m$ for different sample size $N$ }} \\

		\resizebox{!}{!}{
			\begin{tabular}{ p{2.0cm} p{1.0cm}   } 
				\hline
				N  & m   \\
				\hline 
				
				$\leq$ 10 &  2 \\
				11-50     &  6  \\	
				50-100    &  8  \\
				$\geq$ 100 &  10  \\
				\hline
		\end{tabular}}
	\end{center}\label{table9}
	
	\subsection{ Real data application }\label{section7}
	Dataset 1 from Montgomery et al. (2008) has a normal distribution (symmetric model) as a suitable model. \\
	
	Dataset 1: 15.5, 23.75, 8.0, 17.0, 5.5, 19.0, 24.0, 2.5, 7.5, 11.0, 13.0, 3.75, 25.0,9.75, 22.0, 18.0, 6.0, 12.5, 2.0, 21.5.\\
	
	The normal distribution is symmetric. This fact is verified by our test. The value of the test statistics ${\reallywidehat{\Delta}}_{2,2}$ is 0.1531 with an estimated $p$-value 0.2969 when window size $m=2$ and sample size $N=20.$ Our test based on ${\reallywidehat{\Delta}}_{2,2}$ fails in rejecting the null hypothesis even if the significance level is $10\%$. That is because the dataset has a normal distribution as a suitable model.
	
	Dataset 2 from Qiu and Jia (2018b) represents active repair times (in hours) for an airborne communication transceiver. \\
	
	Dataset 2: 0.2, 0.3, 0.5, 0.5, 0.5, 0.5, 0.6, 0.6, 0.7, 0.7, 0.7, 0.8, 0.8, 1.0, 1.0, 1.0, 1.0,1.1, 1.3,1.5,1.5, 1.5, 1.5, 2.0, 2.0, 2.2, 2.5, 3.0, 3.0, 3.3, 3.3, 4.0, 4.0, 4.5, 4.7, 5.0, 5.4, 5.4, 7.0, 7.5, 8.8, 9.0, 10.3, 22.0, 24.5. \\
	
	This data can be fitted by inverse Gaussian (IG) distribution as pointed out in Qiu and Jia  (2018B). IG distribution is not symmetric (see, Xiong et al   (2020) and Qiu and Jia  (2018B)). This fact is verified by our test. The value of the test statistics ${\reallywidehat{\Delta}}_{2,2}$ is 3.6678 with an estimated p-value 0 when window size $m=20$ and sample size $N=45.$ Our test based on ${\reallywidehat{\Delta}}_{2,2}$ succeed
	in rejecting the null hypothesis even if the significance level is small enough, say, $1\%$. This further shows the advantage of our test.\\
	
	Dataset 3: 1.42, 0.84, 2.32, 1.84, 2.4, 0.9, 1.49, 0.87, 1.36, 1.25, 1.25, 1.8, 0.86, 0.04, 0.49, 2.08, 0.58, 0.22, 0.06, 1.7, 2.67, 2.39, 2.32, 2.98, 3.21, 1.99, 1.3, 1.25, 1.76, 1.67, 1.36, 1.57, 1.21, 1.24, 1.62, 0.93, 1.32, 0.86, 1.48, 0.85, 1.23, 1.23, 2.14. \\
	
	Dataset 3 is taken from Sathar and Jose  (2020) and they also used this dataset 3 in testing symmetry. Sathar and Jose (2020) proposed Normal distribution (symmetric model) as a suitable model for this data set. The normal distribution is symmetric. This fact is verified by our test. The value of the test statistics ${\reallywidehat{\Delta}}_{2,2}$ is 0.1545 with an estimated $p$-value 0.2821 when window size $m=3$ and sample size $N=43.$ Our test based on ${\reallywidehat{\Delta}}_{2,2}$ fails in rejecting the null hypothesis even if the significance level is $10\%$. That is because the dataset has a normal distribution as a suitable model.\\
	
	Dataset 4: 99, 61, 86, 113, 96, 99, 83, 57, 80, 79, 75, 70, 15, 62, 87, 95, 81, 71, 44, 13, 52,	97, 146, 52, 52, 29, 108, 135, 102, 48, 66, 90, 22, 72, 176, 107, 84, 83, 37, 67, 83, 36, 49, 39, 102, 66, 154, 72, 63, 83, 77.\\
	
	Dataset 4 is taken from Thomas and Jose (2021). Thomas and Jose  (2021) proposed Burr-type XII distribution (skewed model)	as a suitable model for this data set. Jose and Sathar (2020) also used this dataset 4 in testing symmetry. The value of the test statistics ${\reallywidehat{\Delta}}_{2,2}$ is 6.2144 with an estimated $p$-value 0.0 when window size $m=25$ and sample size $N=51.$  Our test based on ${\reallywidehat{\Delta}}_{2,2}$ succeed
	in rejecting the null hypothesis even if the significance level is small enough, say, $1\%$. \\
	
	Dataset 5 : 0.0518, 0.0518, 0.1009, 0.1009, 0.1917, 0.1917, 0.1917, 0.2336, 0.2336, 0.2336, 0.2733, 0.2733, 0.3467, 0.3805, 0.3805, 0.4126, 0.4431, 0.4719, 0.4719, 0.4993, 0.6162, 0.6550, 0.6550, 0.7059, 0.7211, 0.7356, 0.7623, 0.7863, 0.8178, 0.8810, 0.9337, 0.9404, 0.9732, 0.9858.\\
	
	Dataset 5 transformed vinyl chloride data into uniform
	distribution using probability integral transformation see Xiong et al (2020). The value of the test statistics ${\reallywidehat{\Delta}}_{2,2}$ is 0.0247 with an estimated $p$-value 0.4425 when window size $m=11$ and sample size $N=34.$ Our test based on ${\reallywidehat{\Delta}}_{2,2}$ fails in rejecting the null hypothesis even if the significance level is $10\%$.\\
	
	We considered dataset 6 from Lawless (2011) which represents the quantity of 1000 cycles to failure for electrical appliances in a life test. Gupta and Chaudhary (2023) also used this dataset for testing symmetry.  The value of the test statistics ${\reallywidehat{\Delta}}_{2,2}$ is 0.5776 with an estimated $p$-value 0.0210 when window size $m=2$ and sample size $N=50.$  Our test based on ${\reallywidehat{\Delta}}_{2,2}$ succeeds
	in rejecting the null hypothesis at the significance level is $5\%$. \\
	
	Dataset 6 : 0.014, 0.034, 0.059, 0.061, 0.069, 0.080, 0.123, 0.142, 0.165, 0.210, 0.381, 0.464, 0.479, 0.556, 0.574, 0.839, 0.917, 0.969, 0.991, 1.064, 1.088, 1.091, 1.174, 1.270, 1.275, 1.355, 1.397, 1.477, 1.578, 1.649, 1.702, 1.893, 1.932, 2.001, 2.161, 2.292, 2.326, 2.337, 2.628, 2.785, 2.811, 2.886, 2.993, 3.122, 3.248, 3.715, 3.790, 3.857, 3.912, 4.100.\\

	See Table 11 for the value of test statistics and p-value for different datasets based on the specific window size and sample size of each dataset.\\
	
	\begin{center}
		\noindent{\bf Table 11}.\small{{ Description of models fitted }} 
		
		\resizebox{!}{!}{
			\begin{tabular}{ p{2.0cm} p{1.0cm} p{1.0cm} p{1.0cm} p{2.0cm}  } 
				\hline
				Dataset & N  & m & ${\reallywidehat{\Delta}}_{2,2}$ & p-value  \\
				\hline 
				
				Dataset 1  & 20  & 2  & 0.1531 & 0.2969 \\
				Dataset 2  & 45  & 20 & 3.6678 & 0.0000 \\	
				Dataset 3  & 43  & 3  & 0.1545 & 0.2821 \\
				Dataset 4  & 51  & 25 & 6.2144 & 0.0000 \\
				Dataset 5  & 34  & 11 & 0.0247 & 0.4425  \\
				Dataset 6  & 50  & 2  & 0.5776 & 0.0210 \\
				\hline
		\end{tabular}}
	\end{center}\label{table11}
	If we are testing at a 5\% level of significance then the $p$-value less than 0.05 detects asymmetric nature and a $p$-value more than 0.05 detects the symmetric nature of data.  Table 11 confirms that the newly proposed test determines whether the random sample's distribution is symmetric or asymmetric. The p-values indicate that, at a 5\% level of significance, datasets 2, 4, and 6 do not have symmetry in the distribution of the random sample. Similar to this, a moderate p-value suggests accepting symmetry in the distribution of datasets 1, 3, and 5. As a result, we were able to confirm that the test statistic properly recognised the symmetry in the distribution of the random variable.

	\section{Conclusion and future work}\label{s9conclusion}
	Measures of uncertainty extropy and varextropy are discussed. Estimators for extropy and varextropy are proposed. A test of symmetry is constructed using an estimator of extropy.
	
	One may propose a test of uniformity using Theorem \ref{uniformunique}. Let $X_1, X_2, ..., X_n$ be a random sample from random variable X having support interval [0,1], and $X_{1:n} \leq X_{2:n} \leq \dots \leq X_{n:n}$	are the order statistics of the sample. The hypothesis of interest is $H_0:X$	is uniformly distributed against $H_1:X$ is not uniformly distributed. Let $\Delta_n$ is an estimator of $VJ(X).$ Our proposed estimators converge in probability to $VJ(X),$ that is, our proposed estimators are consistent estimators of $VJ(X)$. Under $H_0,$ $\Delta_n$ converges in probability to 0.  Under $H_1,$ $\Delta_n$  converges in probability to a positive number. Large values of $\Delta_n$ can be regarded as a symptom of non-uniformity and therefore we reject $H_0$ for large values of $\Delta_n.$   Since the test statistic $\Delta_n$ is too complex to determine its distribution under the null hypothesis, One may use Monte Carlo simulation to determine the critical values and power of the test.\\
	\\	
	\noindent \textbf{ \Large Conflict of interest} \\
	\\
	No conflicts of interest are disclosed by the authors.\\
	\\
	\textbf{ \Large Funding} \\
	\\
	Santosh Kumar Chaudhary would like to acknowledge financial support from the Council of Scientific and Industrial Research (CSIR) ( File Number 09/0081(14002) /2022-
	EMR-I ), Government of India. \\

	\section*{Appendix-1}\label{appendix}
	The following steps were used to determine the critical values and compute the power of our proposed test and that of other tests for symmetry at significance level $\alpha=0.10, \ \alpha=0.05, \alpha=0.01:$\\
	(1) we defined a function to calculate the absolute value of ${\reallywidehat{\Delta}}_{2,2}$ .\\
	(2) Generate a sample of size $N$ from the null distribution and compute the test statistics for the sample data;\\
	(3) Repeat Step 2 for 10,000 times and determine the 950th, 975th and 995th quantile respectively of the test statistics as the critical value;\\
	(4) Generate a sample of size $N$ from the alternative distribution and check if the absolute value of the test statistic is greater than the critical value;\\
	(5) Repeat Step 4 for 10,000 times and the percentage of rejection is the power of the test.
	
	\section*{Appendix-2}\label{appendix2}
	The following code in Python is used to calculate the p-value for Dataset 6.\\

	\noindent import numpy as np\\
	def calD2(sample, m, N):\\
	sample.sort()\\
	Junx=-1.0/2/N*sum([(1-2*np.log(1-i/(N+1)))**2*(1-i/(N+1))**4*(sample[min(i+m-1,N-1)]-sample[max(i-m-1,0)])for i in range(1,N+1)])/(2*m/N)\\
	Jlnx=-1.0/2/N*sum([(1-2*np.log(i/(N+1)))**2*(i/(N+1))**4*(sample[min(i+m-1,N-1)] \\ -sample[max(i-m-1,0)])for i in range(1,N+1)])/(2*m/N)\\
	D2=Junx-Jlnx\\
	return D2\\
	\\
	N=50\\
	count2=0\\
	m=2\\
	for i in range(10000):\\
	sample1=np.random.normal(0.0,1.0,N)\\
	sample2=[0.014,0.034,0.059,0.061,0.069,0.080,0.123,0.142,0.165,0.210,0.381,0.464,0.479,0.556, \\ 0.574,0.839,0.917,0.969,0.991,1.064,1.088,1.091,1.174,1.270,1.275,1.355,1.397,1.477,1.578,1.649, \\ 1.702,1.893,1.932,2.001,2.161,2.292,2.326,2.337,2.628,2.785,2.811,2.886,2.993,3.122,3.248,3.715, \\
	3.790,3.857,3.912,4.100]
	if calD2(sample2, m, N) $<$ calD2(sample1, m, N):\\
	count2=count2+1\\
	print("p-value when m=",m,"is",count2/10000)\\
	print("value of test statistics for dataset 6 when m=",m,"is ", calD2(sample2, m, N))\\
	
	\section*{Appendix-3}\label{appendix3}
	The following code in Python is used to calculate the power of the test when an alternative distribution is chi-square with a degree of freedom 3.\\
	
	\noindent import numpy as np \\
	def calD2(sample, m, N):\\
	sample.sort()\\
	Junx=-1.0/2/N*sum([(1-2*np.log(1-i/(N+1)))**2*(1-i/(N+1))**4*(sample[min(i+m-1,N-1)] \\
	-sample[max(i-m-1,0)])for i in range(1,N+1)])/(2*m/N)\\
	Jlnx=-1.0/2/N*sum([(1-2*np.log(i/(N+1)))**2*(i/(N+1))**4*(sample[min(i+m-1,N-1)] \\
	-sample[max(i-m-1,0)])for i in range(1,N+1)])/(2*m/N)\\
	D2=Junx-Jlnx\\
	return D2\\
	\\
	list1=[]\\
	m=2\\
	N=10\\
	for i in range(10000):\\
	sample1=np.random.normal(0.0,1.0,N)\\
	list1.append(calD2(sample1, m, N))\\

	count2=0\\
	for i in range(10000):\\
	sample2=np.random.chisquare(3,N)\\
	if abs(calD2(sample2, m, N))$>$np.quantile(list1, 0.950):\\
	count2=count2+1 \\
	print("power when alpha=0.10",count2/10000)\\	
	\\
	count3=0\\
	for i in range(10000):\\
	sample3=np.random.chisquare(3,N) \\
	if abs(calD2(sample3, m, N))$>$np.quantile(list1, 0.975):\\
	count3=count3+1\\
	print("power when alpha=0.05",count3/10000)\\
	\\
	count4=0\\
	for i in range(10000):\\
	sample4=np.random.chisquare(3,N)\\
	if abs(calD2(sample4, m, N))$>$np.quantile(list1, 0.995):\\
	count4=count4+1\\
	print("power when alpha=0.01",count4/10000)	
	
	\section*{Appendix-4}\label{appendix4}
	The following code in Python is used to calculate critical value when $m=2$ and $N=10.$\\

	\noindent import numpy as np\\
	def calD2(sample, m, N): \\
	sample.sort()\\
	Junx=-1.0/2/N*sum([(1-2*np.log(1-i/(N+1)))**2*(1-i/(N+1))**4*(sample[min(i+m-1,N-1)] \\
	-sample[max(i-m-1,0)])for i in range(1,N+1)])/(2*m/N) \\
	Jlnx=-1.0/2/N*sum([(1-2*np.log(i/(N+1)))**2*(i/(N+1))**4*(sample[min(i+m-1,N-1)]\\
	-sample[max(i-m-1,0)])for i in range(1,N+1)])/(2*m/N)
	D2=Junx-Jlnx \\
	return D2 \\
	
	\noindent list1=[] \\
	m=2 \\
	N=10 \\
	for i in range(10000): \\
	sample1=np.random.normal(0.0,1.0,N) \\
	list1.append(calD2(sample1, m, N)) \\
	print("critical value for 90\% confidence interval means alpha=0.10 is ",np.quantile(list1, 0.950)) \\
	print("critical value for 95\% confidence interval alpha=0.05 is ",np.quantile(list1, 0.975)) \\
	print("critical value for 99\% confidence interval means alpha=0.01 is ",np.quantile(list1, 0.995)) \\

\end{document}